\documentclass[12pt]{article}
\usepackage{amsfonts}
\usepackage{myart}

\input tcilatex
\begin{document}

\title{Logical reloading as overcoming of crisis in geometry}
\author{Yuri A.Rylov}
\date{Institute for Problems in Mechanics, Russian Academy of Sciences,\\
101-1, Vernadskii Ave., Moscow, 119526, Russia.\\
e-mail: rylov@ipmnet.ru\\
Web site: {$http://rsfq1.physics.sunysb.edu/\symbol{126}rylov/yrylov.htm$}\\
or mirror Web site: {$http://gasdyn-ipm.ipmnet.ru/\symbol{126}%
rylov/yrylov.htm$}}
\maketitle

\begin{abstract}
Properties of the logical reloading in the Euclidean geometry are
considered. The logical reloading is a logical operation which replaces one
system of basic concepts of a conception by another system of basic concepts
of the same conception. The logical reloading does not change propositions
of the conception. However, generalizations of the conception are different
for different systems of basic concepts. It is conditioned by the fact, that
some systems of basic concepts contain not only propositions of the
conception, but also some attributes of this conception description.
Properties of the logical reloading are demonstrated in the example of the
proper Euclidean geometry, whose generalization leads to different results
for different system of basic concepts.
\end{abstract}

\section{Introduction}

Geometry studies the shape and mutual disposition of physical bodies,
abstracting from other their properties. After such an abstraction the
physical body turns in a geometrical object, i.e. in some subset of points
of the space. Geometry is a science on a shape and on disposition of
geometrical objects in a space or in a space-time. The space is a set of
points. A geometrical object is a subset of points of the space.

\label{1b}The property of a geometry to be a science on mutual disposition
of geometrical objects in the space or in the space-time will be called
"geometricity". This special term is neccessary, because the contemporary
geometry does not possess the property of "geometricity", in general. In
other words, the contemporary geometry is not always is a science on mutual
disposition \ of geometrical objects. Contemporary mathematics considers a
geometry simply as a logical construction. For instance, the symplectic
geometry is not a science on mutual dispositions of geometric objects. It is
a logical construction, whose form reminds the form of the Euclidean
geometry. In applications of geometry to physics and to mechanics only the
geometricity of the space-time geometry is important. It is of no
importance, whether or not the geometry is a \ logical construction. If the
real space-time geometry is nonaxiomatizable, it means, that it is not a
logical construction. However, such a geometry may not possess the property
of geometricity.\label{1e}

\label{2b}Contemporary mathematicians do not recognize nonaxiomatizable
geometries, which have the property of geometicity, but which are not a
logical construction. This situation should be qualified as a crisis in
geometry \cite{R2005}, which reminds the crisis, when mathematicians did not
recognize non-Euclidean geometries of Lobachevski - Bolyai. \label{2e}

Aforetime the geometry studied disposition of geometrical objects in usual
space. The time was considered as an additional characteristic of the
physical bodies description. After creation of the relativity theory the
space and the time are considered as a united event space (or space-time).
It is a more general approach to description of the event space. Any point
of the event space is an event, which occurs at some place and at some time.

The geometry is described completely, if the distance $\rho $ between any
pair of points belonging to the space is given. The set $\Omega $ of points
with a distance \ $\rho $, given on the set $\Omega $, is known as a metric
space $\mathcal{M}$.

A use of the metric space in the physics and mechanics meets some problems.
These problems lie in the definition of geometric objects in the metric
space $\mathcal{M}$. The distance $\rho $ is supposed to satisfy the
relations%
\begin{equation}
\rho :\Omega \times \Omega \rightarrow \lbrack 0,\infty ),\qquad \rho \left(
Q,P\right) =\rho \left( P,Q\right) ,\qquad \forall P,Q\in \Omega
\label{a1.1}
\end{equation}%
\begin{equation}
\rho \left( P,Q\right) =0,\qquad \text{iff }\ \ \ P=Q  \label{a1.2}
\end{equation}%
\begin{equation}
\rho \left( P,Q\right) +\rho \left( P,R\right) \geq \rho \left( R,Q\right)
,\qquad \forall P,Q,R\in \Omega  \label{a1.3}
\end{equation}%
In the Euclidean geometry the distance has properties (\ref{a1.1}) - (\ref%
{a1.3}).

In the geometry of Minkowski the distance does not possess these properties.
However, it would be very desirable to introduce a metric geometry (or some
analog of metric geometry) for description of the space-time properties,
because the metric geometry is free of such auxiliary concepts as coordinate
system, dimension and such restriction as continuity. Metric geometry
describes the geometric properties in terms of only distance, which is a
true geometric concept.

After removal of the triangle axiom (\ref{a1.3}) the distance geometry
arises \cite{B53}. Blumental failed to construct a straight line in terms of
only distance. He was forced to introduce a straight as a continuous mapping
of interval $\left( 0,1\right) $ onto the space (a point set). Such an
introduction of nonmetric concept of mapping in geometry seems to be
undesirable, because an auxiliary concept is introduced and the distance
geometry ceases to be a pure metric geometry.

In general, a construction of geometrical objects is the main problem of the
metric geometry. One can easily construct sphere and ellipsoid, because in
the Euclidean geometry these geometrical objects are constructed directly in
terms of distance. However, construction of other geometrical objects needs
a use of some auxiliary means. For instance, a definition of a plane
contains a reference to concept of linear independence of vectors. It is not
quite clear, how to introduce this concept in terms of a distance.

A sphere $Sp_{O,P}$ with the center at the point $O$ and the point $P$ on
the surface of the sphere is defined as a set of points $R$%
\begin{equation}
Sp_{O,P}=\left\{ R|\rho \left( O,R\right) =\rho \left( O,P\right) \right\}
\label{a1.4}
\end{equation}%
An ellipsoid $El_{F_{1}F_{2}P}$ with focuses at the points $F_{1},F_{2}$ and
a point $P$ on the surface of the ellipsoid is defined as a set of points $R$%
\begin{equation}
El_{F_{1}F_{2}P}=\left\{ R|\rho \left( F_{1},R\right) +\rho \left(
F_{2},R\right) =\rho \left( F_{1},P\right) +\rho \left( F_{2},P\right)
\right\}  \label{a1.5}
\end{equation}%
If the point $P$ on the surface of ellipsoid coincides with the focus $F_{2}$%
, the ellipsoid $El_{F_{1}F_{2}P}$ degenerates into segment $\mathcal{T}_{%
\left[ F_{1}F_{2}\right] }$ of a straight line.%
\begin{equation}
\mathcal{T}_{\left[ F_{1}F_{2}\right] }=El_{F_{1}F_{2}F_{2}}=\left\{ R|\rho
\left( F_{1},R\right) +\rho \left( F_{2},R\right) =\rho \left(
F_{1},F_{2}\right) \right\}  \label{a1.6}
\end{equation}

In the proper Euclidean geometry the segment $\mathcal{T}_{\left[ F_{1}F_{2}%
\right] }$ has no thickness (it is one-dimensional). However, if the
triangle axiom (\ref{a1.3}) is not satisfied, the set $\mathcal{T}_{\left[
F_{1}F_{2}\right] }$ is a non-one-dimensional surface.

Criterion of one-dimensionality may be formulated in terms of distance. The
section $\mathcal{S}\left( P,\mathcal{T}_{\left[ F_{1}F_{2}\right] }\right) $
is defined as a set of points $R$ 
\begin{equation}
\mathcal{S}\left( P,\mathcal{T}_{\left[ F_{1}F_{2}\right] }\right) =\left\{
R|\rho \left( F_{1},R\right) =\rho \left( F_{1},P\right) \wedge \rho \left(
F_{2},R\right) =\rho \left( F_{2},P\right) \right\} ,\qquad P\in \mathcal{T}%
_{\left[ F_{1}F_{2}\right] }  \label{a1.7}
\end{equation}%
The point $P\in \mathcal{S}\left( P,\mathcal{T}_{\left[ F_{1}F_{2}\right]
}\right) $ in evident way. By definition the segment (\ref{a1.6}) is
one-dimensional (has no thickness), if any section of $\mathcal{T}_{\left[
F_{1}F_{2}\right] }$ consists of one point 
\begin{equation}
\mathcal{S}\left( P,\mathcal{T}_{\left[ F_{1}F_{2}\right] }\right) =\left\{
P\right\} ,\qquad \forall P\in \mathcal{T}_{\left[ F_{1}F_{2}\right] }
\label{a1.8}
\end{equation}%
Let $\mathcal{S}\left( P,\mathcal{T}_{\left[ F_{1}F_{2}\right] }\right) $ be
a section of the segment $\mathcal{T}_{\left[ F_{1}F_{2}\right] }$ at the
point $P\in \mathcal{T}_{\left[ F_{1}F_{2}\right] }$. One can show that the
relation (\ref{a1.8}) takes place, if the distance satisfies the triangle
axiom (\ref{a1.3}). In this case the segment $\mathcal{T}_{\left[ F_{1}F_{2}%
\right] }$ of the straight line can be defined as a line of the shortest
length, and this definition appears to be equivalent to definition (\ref%
{a1.6}).

In other words, if the straight has no thickness, one may use both
definitions of the straight segment $\mathcal{T}_{\left[ F_{1}F_{2}\right] }$%
: (1) the straight line is a shortest line and (2) definition (\ref{a1.6}).
The two definitions are equivalent. However, the definition (\ref{a1.6}) may
be used in the case, when the triangle axiom does not take place, whereas
the first definition becomes to be incorrect.

We have the following problem: "Do such geometries exist, where the triangle
axiom is not satisfied?" Of course, this question is interesting only in
application to the real space-time geometry. It is of no interest for
mathematicians, which may investigate some special class of geometries
(satisfying the triangle axiom), remaining investigation of a more general
geometries for later. Application of metric geometry to the space--time
geometry needs also a refuse from the condition (\ref{a1.2}), which may not
be used in geometries with indefinite metric, for instance, in geometry of
Minkowski.

In the Riemannian space-time geometry we have instead of (\ref{a1.1}) - (\ref%
{a1.3})%
\begin{equation}
\sigma :\Omega \times \Omega \rightarrow \mathbb{R},\qquad \sigma \left(
P,P\right) =0\qquad \sigma \left( Q,P\right) =\sigma \left( P,Q\right)
,\qquad \forall P,Q\in \Omega   \label{a1.9}
\end{equation}%
\begin{equation}
\sqrt{2\sigma \left( P,Q\right) }+\sqrt{2\sigma \left( P,R\right) }\leq 
\sqrt{2\sigma \left( R,Q\right) },\qquad \forall P,Q,R\in \Omega \wedge
\sigma \left( R,Q\right) >0  \label{a1.10}
\end{equation}%
where $\sigma \left( P,Q\right) $ is the world function, connected with the
distance $\rho \left( P,Q\right) $ by means of relation 
\begin{equation}
\sigma \left( P,Q\right) =\frac{1}{2}\rho ^{2}\left( P,Q\right) 
\label{a1.11}
\end{equation}%
In the space-time geometry the world function is always real, and the
distance is positive for timelike interval ($\sigma \left( P,Q\right) >0$),
and it is imaginary for spacelike interval $\left( \sigma \left( P,Q\right)
<0\right) $. In the Riemannian space-time geometry the conditions (\ref{a1.9}%
), (\ref{a1.10}) are fulfilled. Condition (\ref{a1.10}) for spacelike
distances ($\sigma \left( P,Q\right) <0$) is not fulfilled. It is not
important, because spacelike world lines are not used in contemporary
physics.

As soon as the mathematical technique of working with the world function has
been developed \cite{S60,R62,R64}, the question arises: "If the condition (%
\ref{a1.10}) is not fulfilled, is the space-time geometry non-Riemannian, or
is there no geometry at all?" It was a very important question. On one hand,
the metric geometry was insensitive to continuity, or discreteness. It was
insensitive also to dimension of the space-time and to a choice of a
coordinate system. On the other hand, if such non-Riemannian geometries
exist, the habitual axiom of Euclidean (and Riemannian) geometry (the
straight has no thickness) is violated.

We shall not use the term metric geometry with respect to geometry (\ref%
{a1.9}), because the term "metric geometry" is associated with the triangle
axiom imposed on the metric (distance). We shall use the term "physical
geometry" with respect to geometry, which is completely described by the
world function $\sigma $, satisfying the condition (\ref{a1.9}). Another
(more earlier) name of this geometry is "tubular geometry" (T-geometry),
which arose because in T-geometry some straights are substituted by tubes.
There exist such space-time isotropic T-geometries, where timelike tubes
degenerate into one-dimensional straights. For instance, in the space-time
geometry of Minkowski with world function 
\begin{equation}
\sigma _{\mathrm{M}}\left( x,x^{\prime }\right) =\frac{1}{2}g_{ik}\left(
x^{i}-x^{\prime i}\right) \left( x^{k}-x^{\prime k}\right) ,\qquad g_{ik}=%
\text{diag}\left\{ c^{2},-1,-1,-1\right\}  \label{a1.12}
\end{equation}%
the timelike straights are one-dimensional (have no thickness), and motion
of free particles is deterministic.

However, a small deformation (change of $\sigma _{\mathrm{M}}$) of the
space-time of Minkowski transforms timelike straight lines into tubes, and
motion of free particles becomes stochastic. If the space-time deformation
depends on the quantum constant in a proper way, the statistical description
of stochastically moving free particles is equivalent to the quantum
description in terms of the Schr\"{o}dinger equation \cite{R91}. I have
obtained this result only twenty five years ago after the question on
non-Riemannian geometry had appeared.

This fact is connected with the logical reloading, which is a logical
operation. This logical operation realizes a transition from one system of
basic statements of a conception to another system of basic statements of
the same conception. As a logical operation the logical reloading is
essential only at a generalization of the existing conception. Such a
generalizations are rather rear. As far as the logical reloading is used
rather rear, and researchers possess this logical operation rather slightly.

Euclid has created his geometry as a logical construction. The Euclidean
geometry has been taught in such a form for two thousands years. As a result
almost all researchers believe, that a geometry is a logical construction.
But what is connection between space properties and logic? \label{1beg} Is a
geometry a logical construction with necessity? Is the formal logic a
necessary attribute of geometry? Why do mathematicians not recognize
nonaxiomatizable geometries, which do not use the formal logic \cite{R2005}?
This paper is written to answer these questions.

\section{Construction of geometrical objects in Euclidean geometry}

Euclid investigates properties of the space, constructing geometrical
objects. Investigation of the geometrical objects properties meant an
investigation of properties of the space, because one may study geometry,
only via properties of geometrical objects, placed in this space. In other
words, investigation of a set of points is an investigation of properties of
subsets of this point set. Geometry as a logical construction is a
formalization of the process of the geometrical objects construction.

Euclid constructed geometrical objects from blocks. He used three sorts of
blocks: (1) point, (2) segment of straight, (3) angle. Combining these
blocks, Euclid constructed geometrical objects and investigated their
properties. Constructing geometrical objects, Euclid used some rules. Some
part of rules described properties of blocks, another part described process
of combining blocks at construction of geometrical objects. The number of
rules is finite, because the number of block sorts is finite. The Euclidean
space, which had been investigated by Euclid, is uniform and isotropic.
Blocks are not deformed at displacements, and the number of rules,
describing displacement of blocks, is also finite. The rules of working with
blocks admit one to construct geometrical objects from blocks mentally.
Besides, these rules generate the rules of the complicate geometric object
construction from other simpler geometrical objects. The simple geometric
objects, constructed from blocks, are described by the rules known as
axioms. More complicated rules (theorems) of the geometric objects
construction from other geometrical objects are deduced from axioms by rules
of formal logic. As a whole such a construction of geometric objects is
perceived as a logical construction, where rules of the formal logic
reflects rules of construction of geometric objects.

Usually one abstracts from the fact, that the logical construction of a
geometry is a formalization of real construction of geometrical objects from
blocks. The Euclidean geometry is presented directly as a logical
construction. The connection between the geometry and the logical
construction is considered to be so strong, that sometimes a logical
construction, which is not connected with the space properties, is
considered as some kind of a geometry. \label{symplec}For instance, the
symplectic geometry is treated as a kind of a geometry, although it is not
describe space properties. It has only the form of Euclidean geometry with
antisymmetric matrix of metric tensor.

The proper Euclidean geometry is uniform and isotropic. The blocks can be
easily displaced without their deformation, and they are used for
construction of geometrical objects. In the uniform geometry one can use a
finite number of the block sorts. The corresponding logical construction
contains finite number of axioms. All propositions of a uniform geometry can
be deduced from the finite number of axioms, and this geometry may be
qualified as an axiomatizable geometry.

If a geometry $\mathcal{G}$ is not uniform, one cannot use a finite number
of the block sorts, because the blocks are deformed at a displacement. As a
result two similar geometrical objects, constructed in the same way in
different places of the space, will have different properties. We are forced
to use infinite number of block sorts. The number of axioms will be
infinite. Such a geometry should be qualified as nonaxiomatizable geometry.

One should use another way of construction of inhomogeneous geometries. The
inhomogeneous geometry $\mathcal{G}$ is considered as a result of a
deformation of some standard geometry $\mathcal{G}_{\mathrm{st}}$. The
geometry $\mathcal{G}_{\mathrm{st}}$ is axiomatizable, and geometrical
objects in $\mathcal{G}_{\mathrm{st}}$ are constructed from blocks. The
standard geometry $\mathcal{G}_{\mathrm{st}}$ is supposed to be described
completely by the world function $\sigma _{\mathrm{st}}$. In means that all
propositions of $\mathcal{G}_{\mathrm{st}}$ can be expressed in terms of $%
\sigma _{\mathrm{st}}$. In particular, all geometrical objects in $\mathcal{G%
}_{\mathrm{st}}$ can be described in terms of $\sigma _{\mathrm{st}}$.

The inhomogeneous geometry $\mathcal{G}$ and all geometrical objects in $%
\mathcal{G}$ are constructed as result of a deformation of the standard
geometry $\mathcal{G}_{\mathrm{st}}$. It means that in all descriptions of
geometrical objects in $\mathcal{G}_{\mathrm{st}}$ the world function $%
\sigma _{\mathrm{st}}$ is substituted by the world function $\sigma $ of the
geometry $\mathcal{G}$. As a result one obtains descriptions of all
geometrical objects in $\mathcal{G}$ and, hence, one obtains a description
of the geometry $\mathcal{G}$ in terms of the world function $\sigma $.

The method of the geometry $\mathcal{G}$ construction by means of a
deformation of the standard geometry $\mathcal{G}_{\mathrm{st}}$ is called
the deformation principle \cite{R2007}. The proper Euclidean geometry $%
\mathcal{G}_{\mathrm{E}}$ may be used as a standard geometry, because $%
\mathcal{G}_{\mathrm{E}}$ is a axiomatizable physical geometry. The term
"physical geometry" means by definition, that the geometry may be described
completely in terms of the world function $\sigma _{\mathrm{E}}$.

Constructing generalized geometries, one adopts conventionally from Euclid
his method of the geometry construction. However, this method can be used
only for construction of axiomatizable geometries. It does not work at
construction of physical geometries, which are nonaxiomatizable, in general.
We shall use the proper Euclidean geometry $\mathcal{G}_{\mathrm{E}}$ as the
standard geometry $\mathcal{G}_{\mathrm{st}}$, and it means, that we adopt
from Euclid his geometry, but not his method of the geometry construction.

At construction of generalized geometries it is important to compare
geometrical objects in different geometries. In particular, one should be
able to recognize similar objects in different geometries. The geometrical
objects in different physical geometries are considered to be similar, if
they are defined similar, i.e. if they have the same form in terms of the
world functions. For instance, the segment $\mathcal{T}_{\left[ P_{0}P_{1}%
\right] }$ is defined by the relation (\ref{a1.6}) in terms of distance. In
terms of the world function it has the form 
\begin{equation}
\mathcal{T}_{\left[ P_{0}P_{1}\right] }=\left\{ R|\sqrt{2\sigma \left(
P_{0},R\right) }+\sqrt{2\sigma \left( R,P_{1}\right) }=\sqrt{2\sigma \left(
P_{0},P_{1}\right) }\right\}  \label{a2.1}
\end{equation}%
The same form of definition of the segment $\mathcal{T}_{\left[ P_{0}P_{1}%
\right] }$ has in all physical geometries. However, it does not mean, that
properties of the segment $\mathcal{T}_{\left[ P_{0}P_{1}\right] }$ are the
same in all physical geometries, because properties of the world function
are different in different physical geometries. For instance, in the
geometry of Minkowski the timelike segment $\mathcal{T}_{\left[ P_{0}P_{1}%
\right] }$ ($\sigma _{\mathrm{M}}\left( P_{0},P_{1}\right) >0$) is
one-dimensional, i.e. any its section, defined by the relation (\ref{a1.7})
consists of one point 
\begin{equation}
\mathcal{S}\left( P,\mathcal{T}_{\left[ P_{1}P_{2}\right] }\right) =\left\{
P\right\} ,\qquad \forall P\in \mathcal{T}_{\left[ P_{1}P_{2}\right] }
\label{a2.2}
\end{equation}

In the deformed geometry of Minkowski $\sigma _{\mathrm{d}}$, described by
the world function%
\begin{equation}
\sigma _{\mathrm{d}}=\sigma _{\mathrm{M}}+d\cdot \text{sgn}\left( \sigma _{%
\mathrm{M}}\right) ,\qquad d=\frac{\hbar }{2bc}=\text{const}  \label{a2.3}
\end{equation}%
the same timelike segment $\mathcal{T}_{\left[ P_{1}P_{2}\right] }$ is a
three-dimensional surface (tube). In the relation (\ref{a2.3}) $\sigma _{%
\mathrm{M}}$ is the world function of the Minkowski space-time, $\hbar $ is
the quantum constant, $c$ is the speed of the light, and $b$ is some
universal constant. At proper choice of $b$, the world function (\ref{a2.3})
describes the properties of the space-time in microcosm more effective, than
the world function $\sigma _{\mathrm{M}}$. In the space-time geometry $%
\mathcal{G}_{\mathrm{d}}$, described by the world function $\sigma _{\mathrm{%
d}}$, world lines of free microparticles appear to be stochastic.
Statistical description of these stochastic world lines is equivalent to
quantum description in terms of the Schr\"{o}dinger equation \cite{R91}. In
other words, a use of the space-time geometry (\ref{a2.3}) instead of the
geometry of Minkowski admits one to remove quantum principles.

Thus, the proper choice of the space-time geometry admits one to explain
quantum effects as geometrical effects. In this explanation one uses only
principles of classical dynamics. Quantum principles are not introduced, or
they are obtained as corollaries of the physical geometry of the space-time.
At such an approach the number of physical principles reduces. When the
number of basic principles reduces, the physical theory becomes more perfect.

At axiomatic approach to geometry the properties of geometrical objects are
obtained in the form of theorems, deduced from axiomatics. In physical
geometries the properties of geometrical objects are obtained only after
taking into account properties of the world function.

\section{Logical reloading in the proper Euclidean \newline
geometry}

There are three different equivalent representations of the proper Euclidean
geometry \cite{R2007a}: (1) Euclidean representation (E-representation), (2)
vector representation (V-representation) and (3) representation in terms of
world function ($\sigma $-representation). Transformation from one
representation to another one is a logical reloading, when basic concepts of
the representation are changed. The $E$-reprsentation uses three blocks
(point, segment, angle) for construction of geometric objects.

The $V$-representation uses two blocks (point and directed segment, or
vector). Instead of angle one uses additional structure (linear vector
space), which accomplishes function of the angle, describing direction of
vectors.\label{beg2} The angle is constructed from to segments, having a
common point. If one formulates the rules of constructing an angle from two
segments, one may reduce the number of sorts of blocks, remaining only point
and segment. The block "angle" is replaced by the rule of its construction.
As a result one obtains two blocks (point and segment) and some additional
structure, which admits one to construct the angles from two segments. This
additional structure is known as the linear vector space. Reduction of basic
elements (blocks) is a logical reloading, which can be introduced in the
proper Euclidean geometry. This logical reloading from $E$-representation to 
$V$-representation was conditioned by application of Euclidean geometry to
physics and mechanics, where conception of a vector and coordinate system
were used. Concept of the angle was not so essential, because it directivity
of a vector may be described by scalar products of a vector with basic
vectors of the coordinate system. A vector (directed segment) and coordinate
system are attributes of the linear vector space, which is an auxiliary
structure in $V$-representation.\label{end2}

Vector representation is based on the concept of linear vector space, which
contains such concepts as, continuity, dimension, coordinate system, linear
independence. These concepts are necessary for application in physics and
mechanics, where they are used for description of the particle motion and
evolution of the force fields. These concepts are used in $V$-representation
as basic concepts or as properties of the linear vector space.
Conventionally the linear vector space is considered as an attribute of the
Euclidean geometry (but not as an attribute of Euclidean geometry
description).

The $\sigma $-representation contains only one block (point). \label{beg3}
As far as a segment can be constructed of points, it is possible to reduce
the number of block sorts remaining only one sort (point). Elimination of a
segment (vector) is accompanied by the rules of the segment construction
from the points. This reduction of the block sorts leads to the logical
reloading (transition from $V$-representation to $\sigma $-representation).
This transition to $\sigma $-representation is accompanied by introduction
of a new structure (world function $\sigma $), which contains the rules of
the segment construction from points (\ref{a2.1}). The world function
describes a connection of two points of the space. It admits one to
construct all attributes of the linear vector space in terms of the world
function.\label{end3} World function admits one to obtain all concepts of
the $V$-representation (dimension, coordinate system, metric tensor, linear
vector space). Of course, construction of all attributes of the $V$%
-representation is possible, provided the world function is the world
function of the Euclidean space. This world function satisfies some
conditions, which appear to be rather strong.

\begin{definition}
Vector $\mathbf{P}_{0}\mathbf{P}_{1}=\overrightarrow{P_{0}P_{1}}$ is an
ordered set of two points $P_{0},P_{1}$. The point $P_{0}$ is the origin of
the vector and the point $P_{1}$ is its end. The length  of vector $\mathbf{P%
}_{0}\mathbf{P}_{1}$ is
\begin{equation}
\left\vert \mathbf{P}_{0}\mathbf{P}_{1}\right\vert =\sqrt{2\sigma \left(
P_{0},P_{1}\right) }  \label{a3.0}
\end{equation}
\end{definition}

The crucial point of the $\sigma $-representation is the definition of
scalar product in terms of the world function.

\begin{definition}
The scalar product $\left( \mathbf{P}_{0}\mathbf{P}_{1}.\mathbf{Q}_{0}%
\mathbf{Q}_{1}\right) $ of two vectors $\mathbf{P}_{0}\mathbf{P}_{1}$ and $%
\mathbf{Q}_{0}\mathbf{Q}_{1}$ has the form.%
\begin{equation}
\left( \mathbf{P}_{0}\mathbf{P}_{1}.\mathbf{Q}_{0}\mathbf{Q}_{1}\right)
=\sigma \left( P_{0},Q_{1}\right) +\sigma \left( P_{1},Q_{0}\right) -\sigma
\left( P_{0},Q_{0}\right) -\sigma \left( P_{1},Q_{1}\right)   \label{a3.1}
\end{equation}
\end{definition}

If origin of $\mathbf{P}_{0}\mathbf{P}_{1}$ and $\mathbf{Q}_{0}\mathbf{Q}%
_{1} $ is the same $Q_{0}=P_{0}$, the relation (\ref{a3.1}) takes the form%
\begin{equation}
\left( \mathbf{P}_{0}\mathbf{P}_{1}.\mathbf{P}_{0}\mathbf{Q}_{1}\right)
=\sigma \left( P_{0},Q_{1}\right) +\sigma \left( P_{1},P_{0}\right) -\sigma
\left( P_{1},Q_{1}\right)  \label{a3.2}
\end{equation}%
Together with the relation (\ref{a3.0}) the relation (\ref{a3.2}) realizes a
formulation of the cosine theorem. In relations (\ref{a3.0}) - (\ref{a3.2})
the world function $\sigma $ is the world function of the Euclidean geometry.

The necessary and sufficient condition of linear dependence of $n$ vectors ,$%
\mathbf{P}_{0}\mathbf{P}_{1}$, $\mathbf{P}_{0}\mathbf{P}_{2}$,...$\mathbf{P}%
_{0}\mathbf{P}_{n}$, defined by $n+1$ points $\mathcal{P}^{n}\equiv \left\{
P_{0},P_{1},...,P_{n}\right\} $ in the proper Euclidean space, is a
vanishing of the Gram's determinant 
\begin{equation}
F_{n}\left( \mathcal{P}^{n}\right) \equiv \det \left\vert \left\vert \left( 
\mathbf{P}_{0}\mathbf{P}_{i}.\mathbf{P}_{0}\mathbf{P}_{k}\right) \right\vert
\right\vert ,\qquad i,k=1,2,...n  \label{a3.3}
\end{equation}%
Expressing the scalar products $\left( \mathbf{P}_{0}\mathbf{P}_{i}.\mathbf{P%
}_{0}\mathbf{P}_{k}\right) $ in (\ref{a3.3}) via world function $\sigma _{%
\mathrm{E}}$ by means of relation (\ref{a3.2}), we obtain definition of
linear dependence of $n$ vectors $\mathbf{P}_{0}\mathbf{P}_{1}$, $\mathbf{P}%
_{0}\mathbf{P}_{2}$,...$\mathbf{P}_{0}\mathbf{P}_{n}$ in the proper
Euclidean space in the form%
\begin{equation}
F_{n}\left( \mathcal{P}^{n}\right) =0  \label{a3.4}
\end{equation}%
\begin{equation}
F_{n}\left( \mathcal{P}^{n}\right) \equiv \det \left\vert \left\vert \sigma
\left( P_{0},P_{i}\right) +\sigma \left( P_{0},P_{k}\right) -\sigma \left(
P_{i},P_{k}\right) \right\vert \right\vert ,\qquad i,k=1,2,...n  \label{a3.5}
\end{equation}

The necessary and sufficient conditions of the fact, that a physical
geometry, described by the world function $\sigma $, is $n$-dimensional
proper Euclidean geometry, have the form of four conditions.

I. Definition of the dimension of the geometry: 
\begin{equation}
\exists \mathcal{P}^{n}\equiv \left\{ P_{0},P_{1},...P_{n}\right\} \subset
\Omega ,\qquad F_{n}\left( \mathcal{P}^{n}\right) \neq 0,\qquad F_{k}\left( {%
\Omega }^{k+1}\right) =0,\qquad k>n  \label{g2.5}
\end{equation}%
where $F_{n}\left( \mathcal{P}^{n}\right) $\ is the Gram's determinant (\ref%
{a3.5}). Vectors $\mathbf{P}_{0}\mathbf{P}_{i}$, $\;i=1,2,...n$\ are basic
vectors of the rectilinear coordinate system $K_{n}$\ with the origin at the
point $P_{0}$. The metric tensors $g_{ik}\left( \mathcal{P}^{n}\right) $, $%
g^{ik}\left( \mathcal{P}^{n}\right) $, \ $i,k=1,2,...n$\ in $K_{n}$\ are
defined by the relations 
\begin{equation}
\sum\limits_{k=1}^{k=n}g^{ik}\left( \mathcal{P}^{n}\right) g_{lk}\left( 
\mathcal{P}^{n}\right) =\delta _{l}^{i},\qquad g_{il}\left( \mathcal{P}%
^{n}\right) =\left( \mathbf{P}_{0}\mathbf{P}_{i}.\mathbf{P}_{0}\mathbf{P}%
_{l}\right) ,\qquad i,l=1,2,...n  \label{a1.5b}
\end{equation}%
\begin{equation}
F_{n}\left( \mathcal{P}^{n}\right) =\det \left\vert \left\vert g_{ik}\left( 
\mathcal{P}^{n}\right) \right\vert \right\vert \neq 0,\qquad i,k=1,2,...n
\label{g2.6}
\end{equation}

II. Linear structure of the Euclidean space: 
\begin{equation}
\sigma \left( P,Q\right) =\frac{1}{2}\sum\limits_{i,k=1}^{i,k=n}g^{ik}\left( 
\mathcal{P}^{n}\right) \left( x_{i}\left( P\right) -x_{i}\left( Q\right)
\right) \left( x_{k}\left( P\right) -x_{k}\left( Q\right) \right) ,\qquad
\forall P,Q\in \Omega  \label{a1.5a}
\end{equation}%
where coordinates $x_{i}\left( P\right) ,$\ $x_{i}\left( Q\right) ,$ $%
i=1,2,...n$\ of the points $P$ and $Q$\ are covariant coordinates of the
vectors $\mathbf{P}_{0}\mathbf{P}$, $\mathbf{P}_{0}\mathbf{Q}$ respectively,
defined by the relation 
\begin{equation}
x_{i}\left( P\right) =\left( \mathbf{P}_{0}\mathbf{P}_{i}.\mathbf{P}_{0}%
\mathbf{P}\right) ,\qquad i=1,2,...n  \label{b.12}
\end{equation}

III: The metric tensor matrix $g_{lk}\left( \mathcal{P}^{n}\right) $\ has
only positive eigenvalues 
\begin{equation}
g_{k}>0,\qquad k=1,2,...,n  \label{a1.5c}
\end{equation}

IV. The continuity condition: the system of equations 
\begin{equation}
\left( \mathbf{P}_{0}\mathbf{P}_{i}.\mathbf{P}_{0}\mathbf{P}\right)
=y_{i}\in \mathbb{R},\qquad i=1,2,...n  \label{b1.4}
\end{equation}%
considered to be equations for determination of the point $P$\ as a function
of coordinates $y=\left\{ y_{i}\right\} $,\ \ $i=1,2,...n$\ has always one
and only one solution. Conditions I -- IV contain a reference to the
dimension $n$\ of the Euclidean space.

Generalization of the Euclidean geometry in $V$-representation admits one to
consider geometries with indefinite metric tensor (geometry of Minkowski),
or geometry with metric tensor, which is different at different points of
the space. However, in the $V$-representation one can consider only
geometries, which have some dimension, and this dimension is the same at all
points of the space. Besides, in $V$-representation one cannot distinguish
dimension as the number of coordinates, describing manifold from the
dimension as the number of linear independent vectors, although these
concepts are different, in general.

In the $\sigma $-representation one can consider geometries, having no
dimension, or having dimensions, which are different at different points of
the space. This difference between representations arises, because in $V$%
-representation the dimension of space is considered as a primary property
of a geometry, whereas in $\sigma $-representation the dimension is only a
secondary property of a geometry (something like an attribute of the
geometry description in $V$-representation). It is a secondary concept,
determined by the form of the world function.

Logical reloading, transforming $V$-representation of the Euclidean geometry
into the $\sigma $-representation, is very important from the point of view
of possible generalization of the Euclidean geometry. At a generalization of
Euclidean geometry in the $V$-representation the most restrictive properties
(\ref{g2.5}) and (\ref{a1.5a}) of the Euclidean geometry are to be
conserved, because they are properties of the linear vector space, which is
the main structure of the $V$-representation.

At a generalization of Euclidean geometry in $\sigma $-representation the
linear vector space is not used, in general. The equivalence relation
becomes intransitive. Summation of vectors, as well as multiplication of a
vector by a real number become multivariant. In general, a vector cannot be
presented as a sum of its components along the coordinate axes, although
projection $\left( \mathbf{P}_{0}\mathbf{P}_{1}.\mathbf{Q}_{0}\mathbf{Q}%
_{1}\right) /\left\vert \mathbf{Q}_{0}\mathbf{Q}_{1}\right\vert $ of a
vector $\mathbf{P}_{0}\mathbf{P}_{1}$ onto any non-zero vector $\mathbf{Q}%
_{0}\mathbf{Q}_{1}$ is determined uniquely. Although many properties of
vectors in the $\sigma $-representations appear to be multivariant and
unaccustomed, these properties are real properties of space-time geometry.

One should know these real properties of the space-time geometry, because in
the general relativity the space-time geometry is determined by the matter
distribution. One cannot know the space-time geometry previously, and one
should consider all possible geometries. In the general relativity one
supposes, that the space-time geometry can be only a Riemannian geometry.
Thus, the supposition, that the space-time geometry is a Riemannian
geometry, is a mistake from the physical viewpoint. Generalization of the
general relativity on the case of arbitrary physical space-time geometry
shows, that the space-time geometry appears to be non-Riemannian even in the
case of slight gravitational field of a heavy sphere \cite{R2010}.

Mathematics does not consider problems of the geometry application to
physics and to mechanics. Mathematicians may investigate only a part of
possible geometries (for instance, only axiomatizable geometries), and this
consideration of a part of all possible geometries is not a mistake from the
mathematical viewpoint. However, if mathematician believes, that
nonaxiomatizable geometries are impossible and tries to deduce a
nonaxiomatizable geometry from some new axiomatics, it becomes to be a
mistake. The obtained geometry appears to be inconsistent. By the way,
already the Riemannian geometry appears to be inconsistent \cite{R2005}.

\section{Corollaries of the logical reloading}

In application to the Euclidean geometry the logical reloading means a
transition from the conventional $V$-representation to $\sigma $%
-representation. As a result the main structure of the $V$-representation
(linear vector space) is replaced by the structure of the $\sigma $%
-representation (world function). The geometry and all geometrical
quantities are defined via the world function $\sigma $, and only via $%
\sigma $. In particular, vector, which is defined in $V$-representation as
an element of the linear vector space, is defined in $\sigma $%
-representation by relations (\ref{a3.0}) - (\ref{a3.5}). Of course, if
conditions of Euclideaness take place, all results, obtained in $V$%
-representation, coincide with results, obtained in $\sigma $%
-representation. If conditions (\ref{g2.5}) - (\ref{a1.5a}) are not
fulfilled, the geometry ceases to be Euclidean, and a linear vector space
cannot be introduced, in general. However, results, obtained in $\sigma $%
-representation have a sense, and they contain many unexpected properties.

The most important property of the physical geometry, generated by a
violation of conditions (\ref{g2.6}) and (\ref{a1.5a}), is multivariance 
\cite{R2008}. Multivariance of a geometry with respect to vector $\mathbf{P}%
_{0}\mathbf{P}_{1}$ and the point $Q_{0}$ means by definition, that at the
point $Q_{0}$ there are many vectors $\mathbf{Q}_{0}\mathbf{Q}_{1}$, $%
\mathbf{Q}_{0}\mathbf{Q}_{1}^{\prime }$,..., which are equivalent to vector $%
\mathbf{P}_{0}\mathbf{P}_{1}$.

Equivalency $\left( \mathbf{P}_{0}\mathbf{P}_{1}\text{eqv}\mathbf{Q}_{0}%
\mathbf{Q}_{1}\right) $ of vectors $\mathbf{P}_{0}\mathbf{P}_{1}$ and $%
\mathbf{Q}_{0}\mathbf{Q}_{1}$ is defined as follows%
\begin{equation}
\left( \mathbf{P}_{0}\mathbf{P}_{1}\text{eqv}\mathbf{Q}_{0}\mathbf{Q}%
_{1}\right) ,\ \ \text{if }\ \left( \mathbf{P}_{0}\mathbf{P}_{1}.\mathbf{Q}%
_{0}\mathbf{Q}_{1}\right) =\left\vert \mathbf{P}_{0}\mathbf{P}%
_{1}\right\vert \cdot \left\vert \mathbf{Q}_{0}\mathbf{Q}_{1}\right\vert
\wedge \left\vert \mathbf{P}_{0}\mathbf{P}_{1}\right\vert =\left\vert 
\mathbf{Q}_{0}\mathbf{Q}_{1}\right\vert  \label{a4.1}
\end{equation}%
If vector $\mathbf{P}_{0}\mathbf{P}_{1}$ is given, and one looks for vector $%
\mathbf{Q}_{0}\mathbf{Q}_{1}$ at the point $Q_{0}$, which is equivalent to
vector $\mathbf{P}_{0}\mathbf{P}_{1}$, one needs to solve two equations (\ref%
{a4.1}) with respect to position of the point $Q_{1}$. In the proper
Euclidean geometry two equations (\ref{a4.1}) have always a unique solution
independently of dimension of the Euclidean space. It means that at the
point $Q_{0}$ there is one and only one vector $\mathbf{Q}_{0}\mathbf{Q}_{1}$%
, which is equivalent to vector $\mathbf{P}_{0}\mathbf{P}_{1}$. The proper
Euclidean geometry is single-variant. In an arbitrary physical geometry $%
\mathcal{G}$ two equations (\ref{a4.1}) may have many solutions. Then at the
point $Q_{0}$ there are many vectors $\mathbf{Q}_{0}\mathbf{Q}_{1}$, which
are equivalent to vector $\mathbf{P}_{0}\mathbf{P}_{1}$. It means that the
geometry $\mathcal{G}$ is multivariant. In this case the equivalence
relation is intransitive, and geometry $\mathcal{G}$ is nonaxiomatizable.
Thus, nonaxiomatizability of a physical geometry is a corollary of its
multivariance.

The multivariance is a natural property of a physical (nonaxiomatizable)
geometry. Multivariance is absent, when conditions (\ref{g2.6}) and (\ref%
{a1.5a}) are fulfilled. In axiomatizable geometries the property of
multivariance is absent. The axiomatizable geometries have been well studied
and considered usually as "true geometries". However, the axiomatizability
is not a feature of a "true geometry". The axiomatizability is a property of
the most studied geometries. In general, a physical geometry is
multivariant, and the multivariance is a natural property of the space-time
geometry.

As we have seen in introduction, the timelike straight segments in the real
space-time geometry (\ref{a2.3}) are not one-dimensional, and this property
is a corollary of the geometry multivariance with respect to timelike
vectors.

The multivariance leads to a splitting of geometrical objects. We consider
this effect in the example of a circular cylinder. In the proper Euclidean
geometry it is defined by its axis and a point $P$ on its surface. Let $F_{1}
$, $F_{2}$ be two points on the axis of the circular cylinder. The cylinder $%
Cl_{PF_{1}F_{2}}$ is defined as a set of points $R$%
\begin{equation}
Cl_{PF_{1}F_{2}}=\left\{ R|S_{F_{1}F_{2}R}=S_{F_{1}F_{2}P}\right\} 
\label{a4.2}
\end{equation}%
where $S_{F_{1}F_{2}R}$ is the area of the triangle with vertices at the
points $F_{1},F_{2},R$, which is calculated by means of the Heron's formula
via side lengths of the triangle. Let the areas $S_{F_{1}F_{2}R}$ and $%
S_{F_{1}F_{2}P}$ are expressed via world functions of corresponding points.
Let $\mathcal{T}_{\left[ F_{1}F_{2}\right] }$ be the straight segment
between points $F_{1},F_{2}$ and the point $F_{3}\in \mathcal{T}_{\left[
F_{1}F_{2}\right] }$. Let $F_{3}\neq F_{1}$, then in the proper Euclidean
geometry the shape of the circular cylinder depends only on the axis, but
not on a choice of points on this axis, and 
\begin{equation}
Cl_{PF_{1}F_{2}}=Cl_{PF_{1}F_{3}},\qquad F_{3}\in \mathcal{T}_{\left[
F_{1}F_{2}\right] }  \label{a4.3}
\end{equation}%
However, in the multivariant geometry, in general, $Cl_{PF_{1}F_{2}}\neq
Cl_{PF_{1}F_{3}}$ and in the multivariant geometry there are many cylinders,
corresponding to one circular cylinder in the proper Euclidean geometry.
From viewpoint of $V$-representation it is interpreted as a splitting of the
Euclidean cylinder in a multivariant geometry. From viewpoint of $\sigma $%
-representation the fact, that shape of cylinders $Cl_{PF_{1}F_{2}}\ $and $%
Cl_{PF_{1}F_{3}}$ are different, in general, is natural. From this viewpoint
the equation (\ref{a4.3}) means a degeneration of cylinders in the Euclidean
geometry. Interpretation of (\ref{a4.3}) as a degeneration is a more correct
geometrical interpretation, because it does not use such an auxiliary
structure as the linear vector space.

\section{Concluding remarks}

Conventional approach to the space-time geometry, when a geometry is
considered to be a logical construction is poor. To obtain a true
description of the space-time geometry, one needs to realize a logical
reloading and transit to perception of a geometry as a science on shape and
mutual disposition of geometrical objects. Using this approach and
considering the proper Euclidean geometry, one obtains three different
representations of the Euclidean geometry. These representations differ in
the number of block sorts, using for construction of geometrical objects.
Reduction of block sorts is compensated by introduction of an additional
structure, which describes the rules of the eliminated block construction.
In the Euclidean geometry a transition from one representation to another
one may be obtained by means of the formal logic. This circumstance admits
one to interpret the logical reloading as a logical operation.

In the inhomogeneous geometries, where blocks are deformed, a transformation
from one representation to another one becomes to be impossible, in general.
In this case one uses the deformation principle, which works only in $\sigma 
$-representation. After deformation of the standard (Euclidean) geometry the
obtained geometry appears to be nonaxiomatizable, in general.

Logical reloading in the proper Euclidean geometry does not change this
geometry. However, capacities of generalization of the proper Euclidean
geometry are different in different representations of this geometry.
Maximal capacity of generalization appears in the $\sigma $-representation,
where the geometry is described completely in terms of world function $%
\sigma $, which is a function of two points of the space. Generalization of
the proper Euclidean geometry in the $\sigma $-representation admits one to
construct multivariant geometries, which are nonaxiomatizable. Appearance of
multivariant (nonaxiomatizable) geometries shows, that the conventional
approach, when a geometry is a logical construction may be used only in some
special cases. Logical reloading to $\sigma $-representation restores the
old conception of geometry, as a science on a shape and mutual disposition
of geometrical objects. In other words, the logical reloading to the $\sigma 
$-representation is not a new idea. It is a return to old idea of metric
geometry, which was not work without a use of the deformation principle.
Being equipped by the deformation principle, the metric geometry turns into
the physical geometry, which is an excellent tool for description of the
space-time.

A use of the physical geometry admits one to describe space-time geometry of
microcosm, where the geometry is multivariant and cannot be described in
terms of conventional Riemannian geometry. The physical space-time geometry
is effective in cosmology, where the space-time geometry is non-Riemannian
and cannot be described correctly, if one supposes, that the space-time
geometry is Riemannian.

It is interesting, that dynamics of particles is described by finite
difference equations (but not differential). Even dynamic equations for
gravitational field have the form of finite equations (but not
differential). It is connected with the fact, that the physical space-time
geometry may be discrete. Differential equations cannot be used effectively
in the space-time geometry, which may be discrete in some regions.

\end{document}